\documentclass[12pt]{article} 
\usepackage{amsfonts,latexsym,color}
\input{psfig.sty}

\newcounter{environment}[section]
\renewcommand{\theenvironment}{%
\arabic{section}.\arabic{environment}}

{\begin{rm}\refstepcounter{environment}{\textbf\theenvironment\
\bf Definition.~~}}%
{\end{rm}}

{\begin{rm}\refstepcounter{environment}{\textbf\theenvironment\
\bf Example.~~}}%
{\end{rm}}

{\begin{rm}\refstepcounter{environment}{\textbf\theenvironment\
\bf Conjecture.~~}}%
{\end{rm}}

{\begin{rm}\refstepcounter{environment}{\textbf\theenvironment\
\bf Proposition.~~}}%
{\end{rm}}

{\begin{rm}\refstepcounter{environment}{\textbf\theenvironment\
\bf Proposition and Definition.~~}}%
{\end{rm}}

\newenvironment{theorem}%
{\begin{rm}\refstepcounter{environment}{\textbf\theenvironment\
\bf Theorem.~~}}%
{\end{rm}}

{\begin{rm}\refstepcounter{environment}{\textbf\theenvironment\
\bf Corollary.~~}}%
{\end{rm}}

{\begin{rm}\refstepcounter{environment}{\textbf\theenvironment\
\bf Lemma.~~}}%
{\end{rm}}

\def\colon{{:}\;}

\begin{document}
\newcommand{\qbc}[2]{{\left [{#1 \atop #2} \right ]}}
\newcommand{\df}[2]{1-x_{#1}x_{#2}^{-1}}
\newcommand{\pdf}[2]{\left(1-x_{#1}x_{#2}^{-1}\right)}
\newcommand{\beq}{\begin{equation}}
\newcommand{\eeq}{\end{equation}}
\newcommand{\be}{\begin{enumerate}}
\newcommand{\ee}{\end{enumerate}}
\newcommand{\bea}{\begin{eqnarray}}
\newcommand{\eea}{\end{eqnarray}}
\newcommand{\beas}{\begin{eqnarray*}}
\newcommand{\eeas}{\end{eqnarray*}}
\newcommand{\ds}{\displaystyle}
\newcommand{\bm}[1]{{\mbox{\boldmath $#1$}}}
\newcommand{\rr}{\mathbb{R}}
\newcommand{\zz}{\mathbb{Z}}
\newcommand{\st}{\,:\,}
\newcommand{\ps}{\mathbb{Q}[[x]]_0}
\newcommand{\cp}{{\cal P}}
\newcommand{\cf}{{\cal F}}
\newcommand{\gr}{\mathrm{Gr}}
\newcommand{\gl}{\mathrm{GL}}
\newcommand{\sn}{\mathfrak{S}_n}
\newcommand{\cc}{\mathbb{C}}
\newcommand{\qh}{\mathrm{QH}}
\newcommand{\ch}{\mathrm{char}}
\newcommand{\vp}{\varphi}
\newcommand{\dt}{R^{(2)}}
\newcommand{\lm}{\lambda/\mu}
\newcommand{\qtk}{K_{\lambda\mu}(q,t)}
\newcommand{\hs}{\mathrm{Hilb}^n(\cc^2)}
\newcommand{\is}{\mathrm{is}}
\newcommand{\ih}{\mathrm{IH}}
\newcommand{\xp}{X_\cp}
\newcommand{\tm}{\textcolor{magenta}}

\newcommand{\overunder}[2]{
\!\begin{array}{c}
\scriptstyle{#1}\\[-.1in]
-\!\!\!-\!\!\!-\\[-.1in]
\scriptstyle{#2}
\end{array}
\!
}

\newcommand{\wideoverunder}[2]{
\!\begin{array}{c}
\scriptstyle{#1}\\[-.1in]
-\!\!\!-\!\!\!-\!\!\!-\!\!\!-\!\!\!-\!\!\!-\\[-.1in]
\scriptstyle{#2}
\end{array}
\!
}

\newcommand{\overunderarrow}[2]{
\!\begin{array}{c}
\scriptstyle{#1}\\[-.1in]
-\!\!\!-\!\!\!\to\\[-.1in]
\scriptstyle{#2}
\end{array}
\!
}

\newcommand{\wideoverunderarrow}[2]{
\!\begin{array}{c}
\scriptstyle{#1}\\[-.1in]
-\!\!\!-\!\!\!-\!\!\!-\!\!\!-\!\!\!-\!\!\!\to\\[-.1in]
\scriptstyle{#2}
\end{array}
\!
}

\begin{centering}
\textcolor{red}{\Large\bf Recent Developments}\\
\textcolor{red}{\Large\bf in Algebraic Combinatorics}\\[.2in] 
\textcolor{blue}{Richard P. Stanley}\footnote{Partially supported by
  NSF grant 
\#DMS-9988459.}\\ 
Department of Mathematics\\
Massachusetts Institute of Technology\\
Cambridge, MA 02139\\
\emph{e-mail:} rstan@math.mit.edu\\[.1in]
\tm{\small version of 5 February 2004}\\[.2in]
{\textbf Abstract}\\[.1in]
\end{centering}
We survey three recent developments in algebraic combinatorics. The
first is the theory of cluster algebras and the Laurent phenomenon of
Sergey Fomin and Andrei Zelevinsky. The second is the construction of
toric Schur functions and their application to computing three-point
Gromov-Witten invariants, by Alexander Postnikov. The third
development is the construction of intersection cohomology for
nonrational fans by Paul Bressler and Valery Lunts and their
application by Kalle Karu to the toric $h$-vector of a nonrational
polytope. We also briefly discuss the ``half hard Lefschetz theorem''
of Ed Swartz and its application to matroid complexes.

\section{Introduction.}
In a previous paper \cite{rs:ucla} we discussed three recent
developments in algebraic combinatorics. 
In the present paper we
consider three additional topics, namely, the Laurent phenomenon and
its connection with Somos sequences and related sequences, the theory
of toric Schur functions and its connection with the quantum
cohomology of the Grassmannian and 3-point Gromov-Witten invariants,
and the toric $h$-vector of a convex polytope.

\textsc{Note.} The notation $\cc$, $\rr$, and $\zz$,
denotes the sets of complex numbers, real numbers, and integers,
respectively.

\section{The Laurent phenomenon.}
Consider the recurrence
  \beq a_{n-1}a_{n+1}=a_n^2 +(-1)^n,\ n \geq 1, \label{eq:fib} \eeq
with the initial conditions $a_0=0$, $a_1=1$. 
A priori it isn't evident that $a_n$ is an integer for all $n$. 
However, it is easy to check (and is well-known) that $a_n$ is given
by the Fibonacci number $F_n$. 
The recurrence (\ref{eq:fib}) can be ``explained'' by the fact
that $F_n$ is a linear combination of two exponential
functions. 
Equivalently, the recurrence (\ref{eq:fib}) follows from
the addition law for the exponential function $e^x$ or for the sine,
viz.,
$$ 
\sin(x+y) = \sin(x)\cos(y) +\cos(x)\sin(y). 
$$
\indent In the 1980's Michael Somos set out to do something similar
involving the addition law for elliptic functions. 
Around 1982 he discovered a sequence, now known as Somos-6, defined by
quadratic recurrences and seemingly integer valued \cite{somos}.
A number of people generalized Somos-6 to Somos-$N$ for any $N\geq 4$.
The sequences Somos-4 through Somos-7 are defined as follows.
(The definition of Somos-$N$ should then be obvious.)
\beas 
 a_na_{n-4} & = & a_{n-1}a_{n-3} +a_{n-2}^2,\ n\geq 4;\ a_i=1\
     \mathrm{for}\ 0\leq i\leq 3\\
  a_na_{n-5} & = & a_{n-1}a_{n-4} +a_{n-2}a_{n-3},\ n\geq 5;\ a_i=1\
     \mathrm{for}\ 0\leq i\leq 4\\
  a_na_{n-6} & = & a_{n-1}a_{n-5} +a_{n-2}a_{n-4}+a_{n-3}^2,\ n\geq
  6;\\ & & \qquad\qquad\qquad \qquad \qquad\qquad a_i=1\ \mathrm{for}\
  0\leq i\leq 5\\ 
  a_na_{n-7} & = & a_{n-1}a_{n-6} +a_{n-2}a_{n-5}+a_{n-3}a_{n-4},\
  n\geq 7;\\ & & \qquad\qquad\qquad\qquad\qquad\qquad a_i=1\
  \mathrm{for}\ 0\leq i\leq 6.  
\eeas
It was conjectured that all four of these sequences are
integral, i.e., all their terms are integers. 
Surprisingly, however, the terms of Somos-8 are not all integers.
The first nonintegral value is $a_{17}=420514/7$.
Several proofs were quickly given that Somos-4 and Somos-5 are
integral, and independently Hickerson and Stanley showed the
integrality of Somos-6 using extensive computer calculations.
Many other related sequences were either proved or conjectured to be
integral.
For example, Robinson
conjectured that if $1\leq p\leq q\leq r$ and $k=p+q+r$, then the
sequence defined by  
\beq 
 a_na_{n-k} = a_{n-p}a_{n-k+p}+a_{n-q}a_{n-k+q}
     +a_{n-r}a_{n-k+r}, \label{eq:rob} 
\eeq
with initial conditions $a_i=1$ for $0\leq i\leq k-1$, is integral. 
A nice survey of the early history of Somos sequences, including an
elegant proof by Bergman of the integrality of Somos-4 and Somos-5,
was given by Gale \cite{gale}.

A further direction in which Somos sequences can be generalized is the
introduction of parameters. 
The coefficients of the terms of the recurrence can be generic (i.e.,
indeterminates), as first suggested by Gale, and the initial
conditions can be generic.
Thus for instance the generic version of Somos-4 is
\beq 
 a_n a_{n-4} = xa_{n-1}a_{n-3}+ya_{n-2}^2, \label{eq:gen} 
\eeq
with initial conditions $a_0=a$, $a_1=b$, $a_2=c$, and $a_3=d$,
where $x,y,a,b,c,d$ are independent indeterminates. 

Thus $a_n$ is a rational function of the six indeterminates.
A priori the denominator of $a_n$ can be a complicated polynomial, but
it turns out that when $a_n$ is reduced to lowest terms the
denominator is always a monomial, while the numerator is a polynomial
with integer coefficients.
In other words, $a_n\in\zz[x^{\pm 1},y^{\pm 1}, a^{\pm 1},b^{\pm
1},c^{\pm 1},d^{\pm 1}]$, the Laurent polynomial ring over $\zz$ in
the indeterminates $x,y,a,b,c,d$.
This unexpected appearance of Laurent polynomials when more general
rational functions are expected is called by Fomin and Zelevinsky
\cite{lp} the \textbf{Laurent phenomenon}.

Until recently all work related to Somos sequences and the Laurent
phenomenon was of an ad hoc nature. Special cases were proved by
special techniques, and there was no general method for approaching
these problems.
This situation changed with the pioneering work of Fomin and
Zelevinsky \cite{ca}\cite{ca2}\cite{ca3} on \textbf{cluster algebras}.
These are a new class of commutative algebras originally developed in
order to create an algebraic framework for dual-canonical bases and
total positivity in semisimple groups.
A cluster algebra is generated by the union of certain subsets, known
as \textbf{clusters}, of its elements.
Every element $y$ of a cluster is a rational function
$F_y(x_1,\dots,x_n)$ of the elements of any other cluster
$\{x_1,\dots,x_n\}$.
A crucial property of cluster algebras, not at all evident from their
definition, is that $F_y(x_1,\dots,x_n)$ is in fact a \textbf{Laurent
polynomial} in the $x_i$'s.
Fomin and Zelevinsky realized that their proof of this fact could be
modified to apply to a host of combinatorial conjectures and problems
concerning integrality and Laurentness.
Let us note that although cluster algebra techniques have led to
tremendous advances in the understanding of the Laurent phenomena,
they do not appear to be the end of the story.
There are still many conjectures and open problems seemingly not
amenable to cluster algebra techniques.

\begin{figure}
 \centerline{\psfig{figure=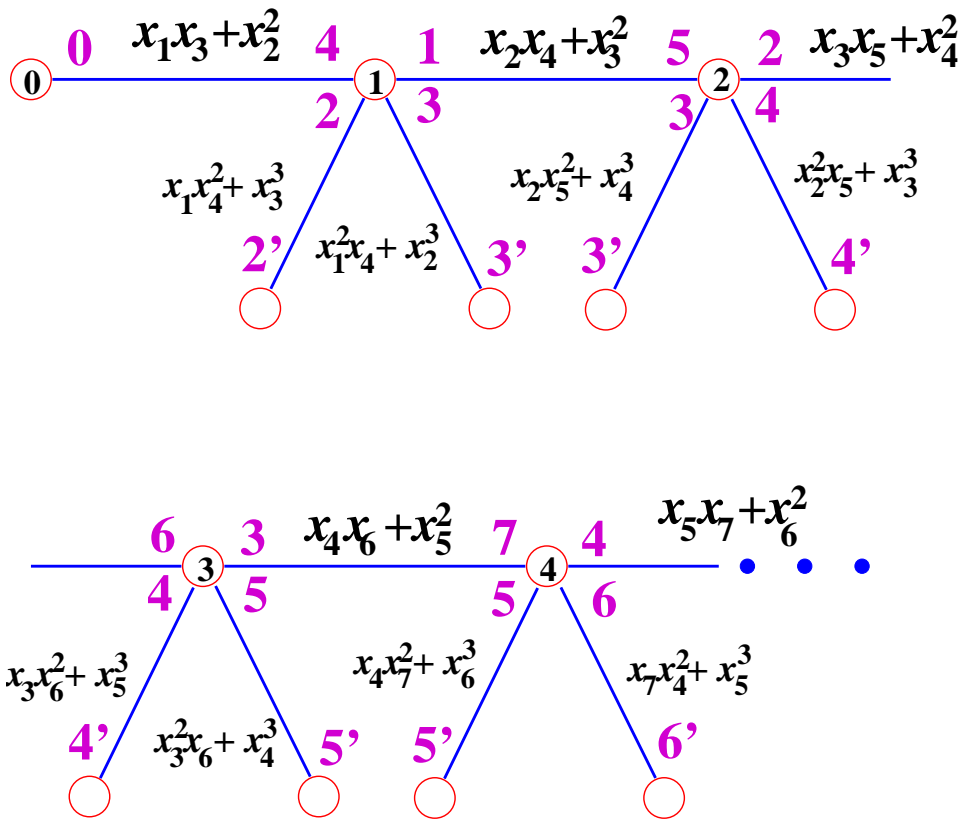}}
\caption{The Somos-4 cluster tree $T$}
\label{fig:somos4}
\end{figure}

We will illustrate the technique of Fomin and Zelevinsky for the
Somos-4 sequence. Consider Figure~\ref{fig:somos4}.
Our variables will consist of $x_0, x_1, x_2, \dots$ and $x_{2'},
x_{3'}, x_{4'}, \dots$.
The figure shows part of an infinite tree $T$, extending to the right.
(We have split the tree into two rows. 
The leftmost edge of the second row is a continuation of the rightmost
edge of the first row.)
The tree consists of a \textbf{spine}, which is an infinite path drawn
at the top, and two legs attached to each vertex of the spine except
the first.
The spine vertices $v=v_i$, $i\geq 0$, are drawn as circles with $i$
inside.
This stands for the set of variables (cluster)
$C_v=\{x_i,x_{i+1},x_{i+2},x_{i+3}\}$.
Each spine edge $e$ has a numerical label $a_e$ on the top left of the
edge, and another $b_e$ on the top right, as well as a polynomial
label $P_e$ above the middle of the edge.
A leg edge $e$ has a numerical label $a_e$ at the top, a polynomial
label $P_e$ in the middle, and a label $b_e=a'_e$ at the bottom.

Moreover, if $e$ is incident to the spine vertex $v$ and leg vertex
$w$, then $w$ has associated with it the cluster
$C_w=(C_v\cup\{x_{a'_e}\})-\{x_{a_e}\}$.
Thus for any edge $e$, if the label $a_e$ is next to vertex $v$ and
the label $b_e$ is next to vertex $w$, then
$C_w=(C_v\cup\{x_{b_e}\})-\{x_{a_e}\}$.

Let $e$ be an edge of $T$ with labels $a_e$, $b_e$, and $P_e$. 
These labels indicate that the variables $x_{a_e}$ and $x_{b_e}$ are
related by the formula
$$ 
x_{a_e}x_{b_e} = P_e. 
$$
(In the situation of cluster algebras, this would be a relation
satisfied by the generators $x_i$.) 
For instance, the leftmost edge of $T$ yields the relation
$$ 
x_0x_4 = x_1x_3+x_2^2. 
$$
In this way all variables $x_i$ and $x'_i$ become rational functions
of the ``initial cluster'' $C_0=\{x_0,x_1,x_2,x_3\}$.

The edge labels of $T$ can be checked to satisfy the following four
conditions:
\begin{itemize}
 \item Every internal vertex $v_i$, $i\geq 1$, has the same
   degree, namely four, and 
   the four edge labels ``next to'' $v_i$ are $i,i+1,i+2,i+3$, the
   indices of the cluster variables associated to $v_i$.
 \item The polynomial $P_e$ does not depend on $x_{a_e}$ and
   $x_{b_e}$, and is not divisible by any variable $x_i$ or $x'_i$.
 \item Write $\bar{P}_e$ for $P_e$ with each variable $x_j$ and $x'_j$
   replaced with $x_{\bar{j}}$, where $\bar{j}$ is the least positive
   residue of $j$ modulo 4. 
   If $e$ and $f$ are consecutive edges of $T$ then the polynomials
   $\bar{P}_e$ and $\bar{P}_{f,0}:= \bar{P}_f|_{x_{\bar{a}_e}=0}$ are
   relatively prime elements of $\zz[x_1,x_2,x_3,x_4]$. 
   For instance, the leftmost two top edges of $T$ yield that
   $x_1x_3+x_2^2$ and $(x_2x_4+x_3^2)|_{x_4=0}=x_3^2$ are coprime.
 \item If $e,f,g$ are three consecutive edges of $T$ such that
   $\bar{a}_e=\bar{a}_g$, then
\beq 
 L\cdot \bar{P}_{f,0}^b\cdot \bar{P}_e =
    \bar{P}_g\mid_{x_{\bar{a}_f} \leftarrow \frac{\bar{P}_{f,0}}{x_{
          \bar{a}_f}}} \label{eq:iji} 
\eeq
  where $L$ is a Laurent monomial, $b\geq 0$, and $x_{\bar{a}_f}
   \leftarrow 
  \frac{\bar{P}_{f,0}}{x_{ \bar{a}_f}}$ denotes the substitution of
  $\frac{\bar{P}_{f,0}}{x_{ \bar{a}_f}}$ for $x_{\bar{a}_f}$.
  For instance, let $e$ be the leftmost leg edge and $f,g$ the second
  and third spine edges. 
  Thus $\bar{a}_e=\bar{a}_g=2$ and
   $\bar{a}_f=1$. Equation~(\ref{eq:iji}) becomes 
$$ 
L\cdot (x_2 x_4+x_3^2)^b_{x_2=0}\cdot (x_1 x_4^2+x_3^3) =
    (x_1x_3+x_4^2)\mid_{x_1\leftarrow\frac{x_3^2}{x_1}}, 
$$
 which holds for $b=0$ and $L=1/x_1$, as desired. 
\end{itemize} 
The above properties may seem rather bizarre, but they are precisely
what is needed to be able to prove by induction that every variable
$x_i$ and $x'_i$ is a Laurent polynomial with integer coefficients in
the initial cluster variables $x_0,x_1,x_2,x_3$ (or indeed in the
variables of any cluster). 
We will not give the proof here, though it is entirely elementary.
Of crucial importance is the periodic nature of the labelled tree $T$.
Each edge is labelled by increasing all indices by one from the
corresponding edge to its left.
This means that the \emph{a priori} infinitely many conditions that
need to be checked are reduced to a (small) finite number.

It follows from the relations $x_ix_{i+4}=x_{i+1}x_{i+3}+x_{i+2}^2$
that $x_n$ is just the $n$th term of Somos-4 with the generic initial
conditions $x_0,x_1,x_2,x_3$. 
Since $x_n$ is a Laurent polynomial with integer coefficients in the
variables $x_0,x_1,x_2,x_3$, if we set $x_0=x_1=x_2=x_3=1$ then $x_n$
becomes an integer.
In this way the integrality of the original Somos-4 sequence
is proved by Fomin and Zelevinsky. 

By similar arguments Fomin and Zelevinsky prove a host of other
integrality theorems, as mentioned above. 
In particular, they prove the integrality of Somos-5, Somos-6,
Somos-7, and the Robinson recurrence (\ref{eq:rob}) by this method.
This gives the first proof of the integrality of Somos-7 (and the
first published proof for Somos-6), as well as a proof of Robinson's
conjecture.
By a refinement of the argument of Fomin and Zelevinsky for Somos-4,
David Speyer \cite{speyer} has shown that for the generic Somos-4
sequence (\ref{eq:gen}) (with generic initial conditions), and
similarly for generic Somos-5, the coefficients of the Laurent
polynomial $x_n$ are \emph{nonnegative}.
Nonnegativity remains open for generic Somos-6 and Somos-7.
The reader might find it instructive to modify (straightforwardly) the
graph $T$ of Figure~\ref{fig:somos4} to prove the following
\cite[Example~3.3]{lp}.

\begin{theorem} \label{thm:gens4}
\emph{Let $a,b$, and $c$ be positive integers, and let the sequence
  $y_0,y_1,\dots$ satisfy the recurrence}
  $$ y_k = \frac{y_{k-3}^ay_{k-1}^c+y_{k-2}^b}{y_{k-4}}. $$
\emph{Then each $y_i$ is a Laurent polynomial with integer
  coefficients in the initial terms $y_0,y_1,y_2,y_3$.}
\end{theorem}

Once the integrality of a recurrence is proved, it is natural to ask
for a combinatorial proof. 

In the case of Somos-4, we would like to
give a combinatorial interpretation to the terms $a_n$ and from this a
combinatorial proof of the recurrence
$$ 
a_n a_{n-4} = a_{n-1}a_{n-3} + a_{n-2}^2.
$$
A clue as to how this might be done comes from the observation that
the rate of growth of $a_n$ is roughly quadratically
exponential. 
Indeed, the function $\alpha^{n^2}$ satisfies the Somos-4 recurrence
if $\alpha^8=\alpha^2+1$.
A previously known enumeration problem whose solution grows
quadratically exponentially arises from the theory of matchings or
domino (dimer) tilings.
Let $G$ be a finite graph, which we assume for convenience has no
loops (vertices connected to themselves by an edge).
A \textbf{complete matching} of $G$ consists of a set of
vertex-disjoint edges that cover all the vertices.
Thus $G$ must have an even number $2m$ of vertices, and each complete
matching contains $m$ edges.

\begin{figure}\centerline{\psfig{figure=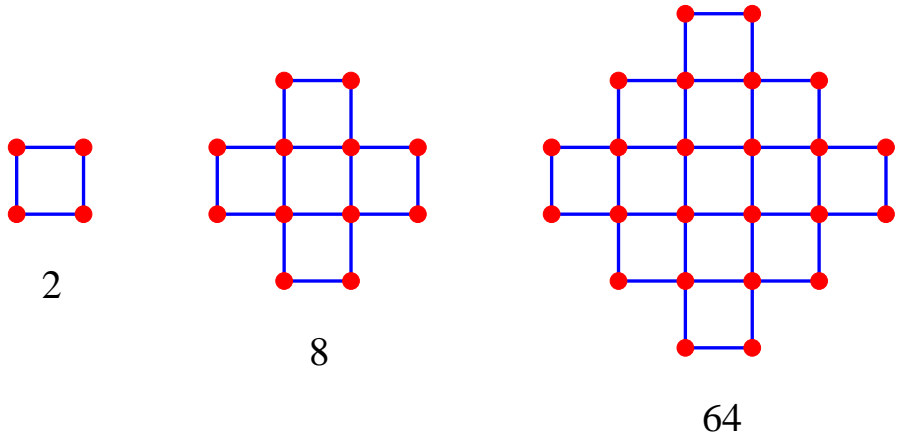}}
\caption{The Aztec diamond graphs AZ$_n$ for $1\leq n\leq 3$}
\label{fig:aztec}
\end{figure}

\begin{figure}
\centerline{\psfig{figure=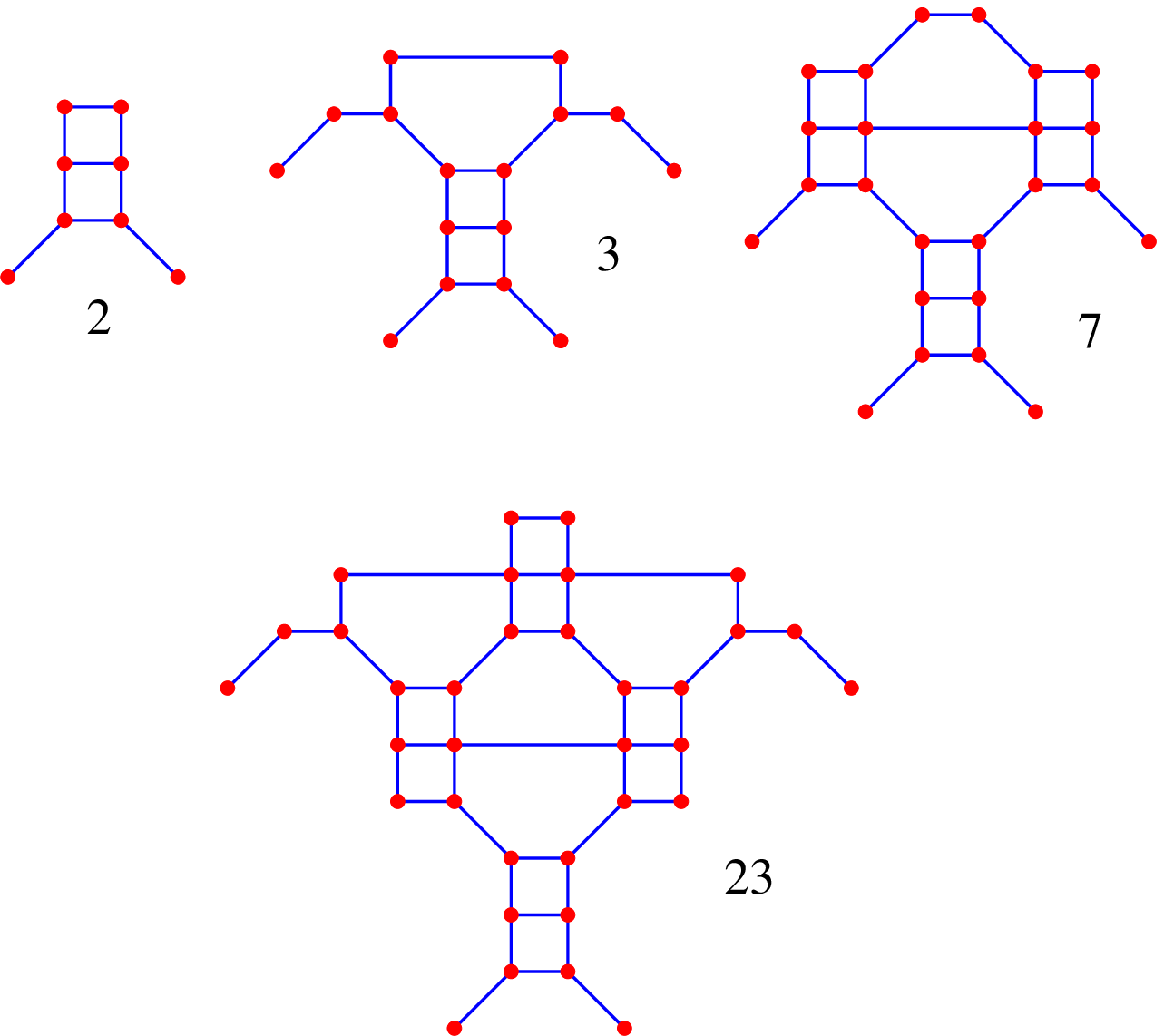}}
\caption{The Somos-4 graphs}
\label{fig:somos}
\end{figure}

Figure~\ref{fig:aztec} shows a sequence of graphs AZ$_1$, AZ$_2$,
AZ$_3, \dots$, whose general definition should be
clear from the figure. 
These graphs were introduced by Elkies, Kuperberg, Larsen, and Propp
\cite{eklp}\cite{eklp2}, who called them (essentially) \textbf{Aztec
diamond graphs}.
They give four proofs that the number of complete matchings of AZ$_n$
is $2^{{n+1\choose 2}}$.
Since this number grows quadratically exponentially, Jim Propp got the
idea that the terms $a_n$ of Somos-4 might count the number of
complete matchings in a planar graph $S_n$ for which the Somos-4
recurrence could be proved combinatorially.
The undergraduate research team REACH \cite{reach}, directed by Propp,
and independently Bousquet-M\'elou, Propp, and West \cite{b-p-w}
succeeded in finding such graphs $S_n$ in the spring of 2002
\cite{reach}. 
Figure~\ref{fig:somos} shows the ``Somos-4 graphs'' $S_4,S_5,S_6,S_7$
along with their number of complete matchings.

\section{Gromov-Witten invariants and toric\\ Schur functions.}
\label{sec:toric} 
Let $\gr_{kn}$ denote the set of all $k$-dimensional subspaces of the
$n$-dimensional complex vector space $\cc^n$. 
We call $\gr_{kn}$ the \textbf{Grassmann variety} or
\textbf{Grassmannian}. It has the structure of a complex projective
variety of dimension $k(n-k)$ and is naturally embedded in complex
projective space $P^{{n\choose k}-1}(\cc)$ of dimension ${n\choose
k}-1$.
The cohomology ring $H^\ast(\gr_{kn})=H^\ast(\gr_{kn};\zz)$ is the
fundamental object for the development of classical \textbf{Schubert
calculus}, which is concerned, at the enumerative level, with counting
the number of linear subspaces that satisfy certain geometric
conditions.
For an introduction to Schubert calculus see
\cite{fulton}\cite{kleiman}, and for connections with combinatorics
see \cite{rs:sc}.
In this section we explain some recent results of Alexander Postnikov
\cite{post} on a quantum deformation of $H^\ast(\gr_{kn})$.
Further details and references may be found in \cite{post}.

A basis for the cohomology ring $H^\ast(\gr_{kn})$ consists of
\textbf{Schubert classes} $\sigma_\lambda$, where $\lambda$ ranges
over all partitions whose shape fits in a $k\times (n-k)$ rectangle,
i.e, $\lambda=(\lambda_1,\dots,\lambda_k)$ where $n-k\geq
\lambda_1\geq \cdots \geq \lambda_k\geq 0$. 
Let $P_{kn}$ denote the set of all such partitions, so
$$
\#P_{kn}=\mathrm{rank}\, H^\ast(\gr_{kn}) = {n\choose k}. 
$$
The Schubert classes $\sigma_\lambda$ are the cohomology classes of the
\textbf{Schubert varieties} $\Omega_\lambda\subset \gr_{kn}$, which are
defined by simple geometric conditions, viz., certain bounds on the
dimensions of the intersections of a subspace $X\in\gr_{kn}$ with the
subspaces $V_i$ in a fixed flag $\{0\}=V_0\subset
V_1\subset\cdots\subset V_n=\cc^n$. 
Multiplication in the ring $H^\ast(\gr_{kn})$ is given by
\beq 
 \sigma_\mu \sigma_\nu
    =\sum_{\lambda\in P_{kn}} c_{\mu\nu}^\lambda \sigma_\lambda,
    \label{eq:grmul} 
\eeq 
where $c_{\mu\nu}^\lambda$ is a \textbf{Littlewood-Richardson
coefficient}, described combinatorially by the famous
\textbf{Littlewood-Richardson rule} (e.g., \cite[Appendix~A]{ec2}).
Thus $c_{\mu\nu}^\lambda$ has a geometric interpretation as the
intersection number of the Schubert varieties
$\Omega_\mu,\Omega_\nu,\Omega_\lambda$.
More concretely, 
\beq
 c_{\mu\nu}^\lambda = \#\left(\tilde{\Omega}_\mu\cap
  \tilde{\Omega}_\nu\cap \tilde{\Omega}_{\lambda^\vee}\right),
 \label{eq:geomlr} 
\eeq 
the number of points of $\gr_{kn}$ contained in the intersection
$\tilde{\Omega}_\mu\cap\tilde{\Omega}_\nu\cap\tilde{\Omega}_{\lambda^\vee}$,
where $\tilde{\Omega}_\sigma$ denotes a generic translation of
$\Omega_\sigma$ and $\lambda^\vee$ is the \textbf{complementary
partition} $(n-k-\lambda_k,\dots,n-k-\lambda_1)$.
Equivalently, $c_{\mu\nu}^\lambda$ is the number of $k$-dimensional
subspaces of $\cc^n$ satisfying all of the geometric conditions
defining $\tilde{\Omega}_{\lambda^\vee}$, $\tilde{\Omega}_\mu$, and
$\tilde{\Omega}_\nu$.

The cohomology ring $H^\ast(\gr_{kn})$ can be deformed into a
``quantum cohomology ring'' $\qh^\ast(\gr_{kn})$, which specializes to
$H^\ast(\gr_{kn})$ by setting $q=0$. 
More precisely, let $\Lambda_k$ denote the ring of symmetric
polynomials over $\zz$ in the variables $x_1,\dots,x_k$. Thus
$$ 
\Lambda_k = \zz[e_1,\dots,e_k], 
$$
where $e_i$ is the $i$th elementary symmetric function in the
variables $x_1,\dots, x_k$, viz,
$$ 
e_i = \sum_{1\leq j_1<j_2<\cdots<j_i\leq k}x_{j_1}x_{j_2}\cdots\
     x_{j_k}. 
$$
Then we have the ring isomorphism
\beq H^\ast(\gr_{kn})\cong \Lambda_k/(h_{n-k+1},\dots,h_n), 
     \label{eq:grcoho} 
\eeq
where $h_i$ denotes a complete homogeneous symmetric function (the sum
of all distinct monomials of degree $i$ in the variables
$x_1,\dots,x_k$).  
The isomorphism (\ref{eq:grcoho}) associates the Schubert class
$\sigma_\lambda\in H^*(\gr_{kn})$ with the (image of) the
\textbf{Schur function} $s_\lambda(x_1,\dots,x_k)$.
The Schur functions $s_\lambda$, where $\lambda$ has at most $k$
parts, form a basis for the abelian group $\Lambda_k$ with many
remarkable combinatorial and algebraic properties
\cite{macd}\cite[Ch.~7]{ec2}.
An important property of Schur functions is their \textbf{stability},
viz., $s_\lambda(x_1,\dots,x_k,0)=s_\lambda(x_1,\dots,x_k)$, which
allows us to define the Schur function
$s_\lambda=s_\lambda(x_1,x_2,\dots)$ in infinitely many variable $x_i$
by $s_\lambda=\lim_{k\rightarrow\infty} s_\lambda(x_1,\dots,x_k)$.
It follows from (\ref{eq:grmul}) and (\ref{eq:grcoho}) that
\beq 
 s_\mu s_\nu=\sum_\lambda c_{\mu\nu}^\lambda s_\lambda, 
   \label{eq:lrdef} 
\eeq
the ``usual'' definition of the Littlewood-Richardson coefficients
$c_{\mu\nu}^\lambda$. 

The quantum cohomology ring $\qh^\ast(\gr_{kn})$ differs from
$H^\ast(\gr_{kn})$ in just one relation: we must enlarge the
coefficient ring to $\zz[q]$ and replace the relation $h_n=0$ with
$h_n=(-1)^{k-1}q$. 
Thus
\beq 
 \qh^\ast(\gr_{kn}) \cong \Lambda_k\otimes \zz[q]/
    (h_{n-k+1},\dots,h_{n-1},h_n+(-1)^kq). \label{eq:qgrcoho} 
\eeq
A basis for $\qh^\ast(\gr_{kn})$ remains those $\sigma_\lambda$ whose
shape fits in a $k\times (n-k)$ rectangle, and under the isomorphism
(\ref{eq:qgrcoho}) $\sigma_\lambda$ continues to correspond to the
Schur function $s_\lambda$. 
Now, however, the usual multiplication $\sigma_\mu \sigma_\nu$ of
Schubert classes has been deformed into a ``quantum multiplication''
$\sigma_\mu\ast \sigma_\nu$.
It has the form 
\beq 
 \sigma_\mu\ast\sigma_\nu = \sum_{d\geq 0} \sum_{{\lambda\vdash
      |\mu| +|\nu| -dn\atop\lambda\in P_{kn}}}q^d
    C_{\mu\nu}^{\lambda,d} \sigma_\lambda, \label{eq:gw} 
\eeq 
where $C_{\mu\nu}^{\lambda,d}\in\zz$. 
The geometric significance of the coefficients
$C_{\mu\nu}^{\lambda,d}$ (and the motivation for defining
$\qh^\ast(\gr_{kn})$ in the first place) is that they count the number
of rational curves of degree $d$ in $\gr_{kn}$ that meet fixed generic
translates of the Schubert varieties $\Omega_{\lambda^\vee}$,
$\Omega_\mu$, and $\Omega_\nu$.
(Naively, a \textbf{rational curve} of degree $d$ in $\gr_{kn}$ is a
set $C=\{(f_1(s,t), f_2(s,t), \dots,f_{{n\choose k}}(s,t))\in
P^{{n\choose k}-1}(\cc)\colon s,t\in \cc\}$, where $f_1(x,y),\dots,
f_{{n\choose k}}(x,y)$ are homogeneous polynomials of degree $d$ such
that $C\subset \gr_{kn}$.)
Since a rational curve of degree 0 in $\gr_{kn}$ is just a point of
$\gr_{kn}$ we recover in the case $d=0$ the geometric interpretation
(\ref{eq:geomlr}) of ordinary Littlewood-Richardson coefficients
$c_{\mu\nu}^\lambda = C_{\mu\nu}^{\lambda,0}$.
The numbers $C_{\mu\nu}^{\lambda,d}$ are known as (3-point)
\textbf{Gromov-Witten invariants}.
From this geometric interpretation of Gromov-Witten invariants it
follows that $C_{\mu\nu}^{\lambda,d}\geq 0$.
No algebraic or combinatorial proof of this inequality is known using
equation (\ref{eq:qgrcoho}) as the definition of $\qh^\ast(\gr_{kn})$,
and it is a fundamental open problem to find such a proof.

The primary contribution of Postnikov is a combinatorial description
of a new generalization of (skew) Schur functions whose expansion
into Schur functions has coefficients which are the Gromov-Witten
invariants. 
This description leads to a better understanding of earlier results as
well as a host of new results.

We begin by reviewing the combinatorial definition of skew Schur
functions. 

Let $\mu,\lambda$ be partitions with $\mu\subseteq
\lambda$, i.e., $\mu_i\leq \lambda_i$ for all $i$. 
The pair $(\mu,\lambda)$ is called a \textbf{skew partition}, often
denoted $\lm$.
The \textbf{diagram} of $\lm$ consists of the diagram of $\lambda$
with $\mu$ removed (from the upper-left corner).
For example, the diagram of $(4,4,3,1)/(2,1,1)$ is given by 

\vspace{.2in}
\centerline{\psfig{figure=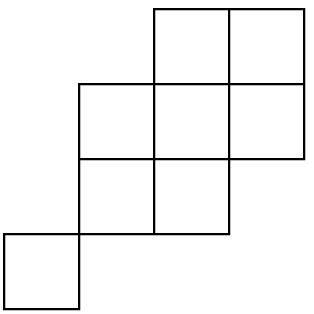}}

\noindent We sometimes call the diagram of $\lm$ a \textbf{skew shape}. 

A \textbf{semistandard Young tableau} (SSYT) of shape $\lm$
consists of the diagram of $\lm$ with a positive integer
placed in each square, so that the rows are weakly increasing and
columns strictly increasing. 
If $T$ is an SSYT with $a_i$ occurences of the entry $i$, then write
\beq 
 x^T=x_1^{a_1}x_2^{a_2}\cdots. \label{eq:xt} 
\eeq
Hence the total degree of the monomial $x^T$ is
$|\lm|=|\lambda|-|\mu|$. Define the \textbf{skew Schur function}
$s_{\lm}$ by
$$ 
s_{\lm} = \sum_T x^T, 
$$
where $T$ ranges over all SSYT of shape $\lm$. 
The basic facts concerning $s_{\lm}$ are the following:
  \begin{itemize}
   \item Let $\mu=\emptyset$, so $\lm$ is just the ``ordinary''
     partition $\lambda$. 
     Then $s_{\lambda/\emptyset}=s_\lambda$, the ``ordinary'' Schur
     function. 
  \item $s_{\lm}$ is a symmetric function, whose expansion in terms of
    Schur functions is given by 
\beq  
 s_{\lm}=\sum_\nu c_{\mu\nu}^\lambda s_\nu, \label{eq:skews} 
\eeq
where $c_{\mu\nu}^\lambda$ denotes a Littlewood-Richardson coefficient
(see e.g.\ \cite[(7.164)]{ec2}).
  \end{itemize}

\begin{figure}
 \centerline{\psfig{figure=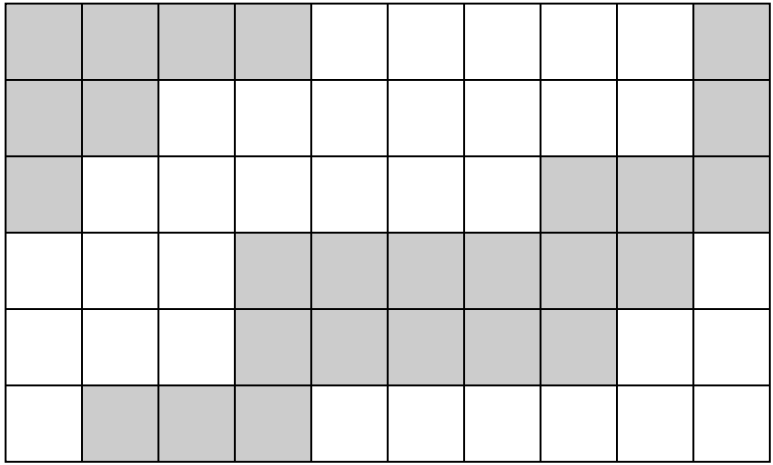}}
\caption{A toric shape in a $6\times 10$ rectangle}
\label{fig:toricsh}
\end{figure}

\begin{figure}
 \centerline{\psfig{figure=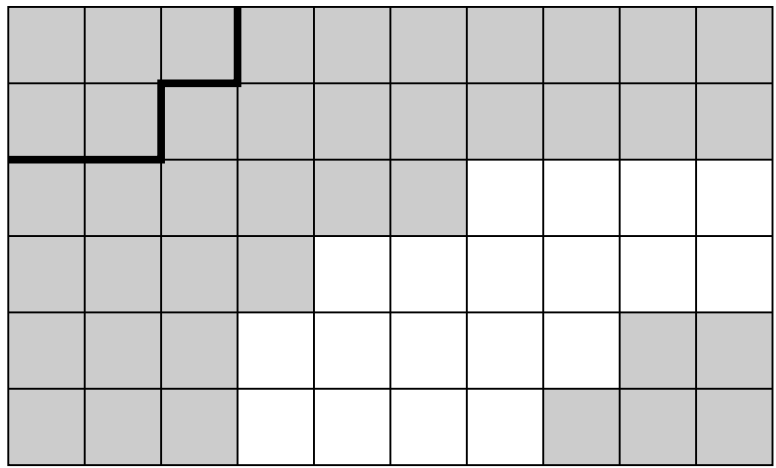}}
\caption{A toric shape with loops}
\label{fig:torloop}
\end{figure}

We see from equations (\ref{eq:lrdef}) and (\ref{eq:skews}) that
$c_{\mu\nu}^\lambda$ has two ``adjoint'' descriptions, one as a
coefficient in a product of Schur functions and one as a coefficient
in a skew Schur function. 
We have already seen in equation~(\ref{eq:gw}) the quantum analogue of
(\ref{eq:lrdef}) so now we would like to do the same for
(\ref{eq:skews}).
We do this by generalizing the definition of a skew shape to a
\textbf{toric shape}.
(Postnikov develops the theory of toric shapes within the more general
framework of \textbf{cylindric shapes}, but we will deal directly with
toric shapes.)
An ordinary skew shape is a certain subset of squares of a $k\times
(n-k)$ rectangle.
A toric shape is a certain subset $\tau$ of squares of a $k\times
(n-k)$ \emph{torus} (which we regard as a $k\times (n-k)$ rectangle
with the left and right edges identified, and the top and bottom edges
identified).
Namely, (a) each row and column of $\tau$ is an unbroken line of
squares (on the torus), and (b) if $(i,j),(i+1,j+1)\in\tau$ (taking
indices modulo $(k,n-k)$), then $(i+1,j),(i,j+2)\in\tau$.
Figure~\ref{fig:toricsh} illustrates a typical toric shape (taken from
\cite{post}).
(If a row or column forms a loop around the torus, then we must also
specify which square is the initial square of the row or column, and
these specifications must satisfy a natural consistency condition
which we hope Figure~\ref{fig:torloop} makes clear.
Henceforth we will the ignore the minor modifications needed in our
definitions and results when the shape contains toroidal loops.)

\begin{figure}
 \centerline{\psfig{figure=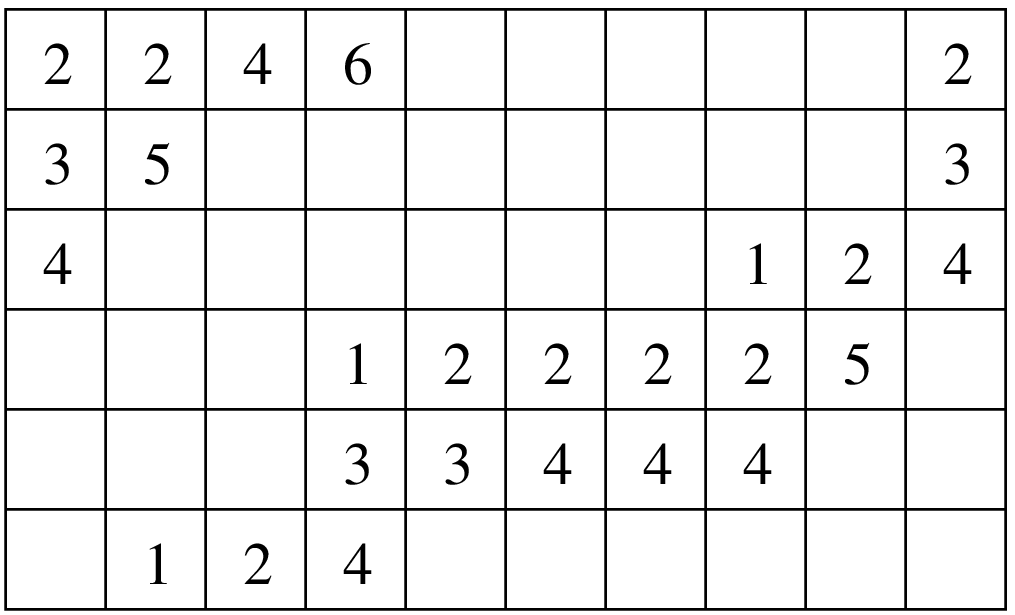}}
\caption{A semistandard toric tableau}
\label{fig:sstt}
\end{figure}

If $\tau$ is a toric shape, then we define a \textbf{semistandard toric
tableau} (SSTT) of shape $\tau$ in exact analogy to the definition
of an SSYT: place positive integers into the squares of $\tau$ so
that every row is weakly increasing and every column strictly
increasing. 
Figure~\ref{fig:sstt} shows an SSTT of the shape $\tau$ given by
Figure~\ref{fig:toricsh}.

\begin{figure}
 \centerline{\psfig{figure=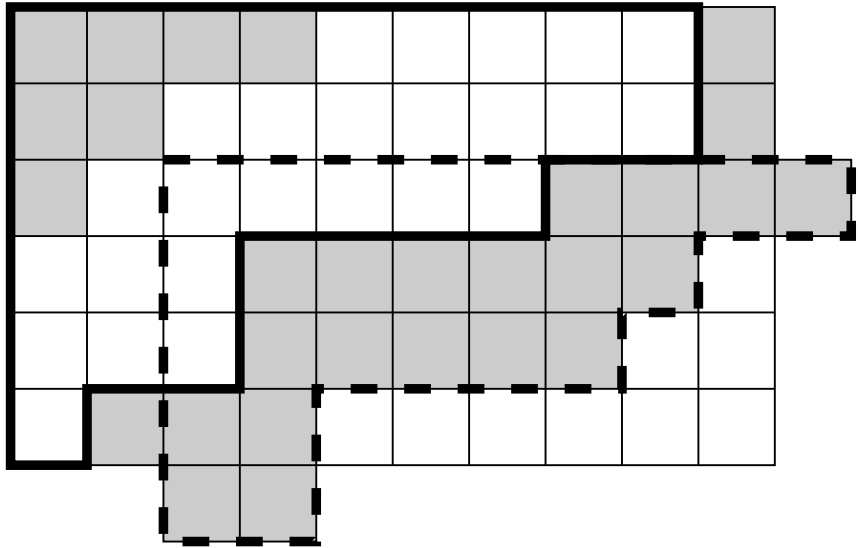}}
\caption{The toric shape $(9,7,6,2,2,0)/2/(9,9,7,3,3,1)$}
\label{fig:toriclm}
\end{figure}

We now explain a method of indexing a toric shape $\tau$ by a triple
$\lambda/d/\mu$, where $\lambda$ and
$\mu$ are partitions and $d\geq 0$. 
We will illustrate this indexing with the toric shape of
Figure~\ref{fig:toricsh}.
In Figure~\ref{fig:toriclm} we have placed the ordinary shape
$\mu=(9,9,7,3,3,1)$, outlined with dark solid lines, on the $6\times
10$ torus $R$.
It is the largest shape contained in $R$ whose intersection with
$\tau$ is an ordinary shape.
Similarly, translated $d=2$ diagonal steps from $\mu$ is the shape
$\lambda= (9,7,6,2,2,0)$, outlined with dark broken lines (and drawn
for clarity to extend beyond $R$ but regarded as being on the torus
$R$).
It is the largest shape whose upper-left hand corner is a diagonal
translation of the upper left-hand corner of $\mu$ and whose
intersection with the complement of $\tau$ is a subshape of $\lambda$.
Thus we rewrite the shape $\tau$ as  
$$ 
\tau = \lambda/d/\mu = (9,7,6,2,2,0)/2/(9,9,7,3,3,1). 
$$
This representation is not unique since any square of $R$ could be
taken as the upper-right corner, but this is irrelevant to the
statement of Theorem~\ref{thm:post} below.

We now define the \textbf{toric Schur function} $s_{\lambda/d/\mu}$
exactly in analogy to the skew Schur function $s_{\lambda/\mu} =
s_{\lambda/0/\mu}$, viz.,
$$ 
s_{\lambda/d/\mu} =\sum_T x^T, 
$$
summed over all SSTT of shape $\lambda/d/\mu$, where $x^T$ is defined
exactly as in (\ref{eq:xt}). 
The remarkable main theorem of Postnikov \cite[Thm.\ 6.3]{post} is the
following.

\begin{theorem} \label{thm:post}
\emph{Let $\lambda/d/\mu$ be a toric shape contained in a $k\times
(n-k)$ torus. Then}
\beq s_{\lambda/d/\mu}(x_1,\dots,x_k) = \sum_{\nu\in P_{kn}}
 C_{\mu\nu}^{\lambda,d}s_\nu(x_1,\dots,x_k). \label{eq:gwmain}
\eeq
\end{theorem}

\textsc{Note.} The above theorem shows in particular that the
expansion of the toric Schur function
$s_{\lambda/d/\mu}(x_1,\dots,x_k)$ into Schur functions 
has nonnegative coefficients, i.e., $s_{\lambda/d/\mu}(x_1,\dots,x_k)$
is \textbf{Schur positive}. 
We mentioned above that no ``direct'' proof using (\ref{eq:gw}) was
known that $C_{\mu\nu}^{\lambda,d}\geq 0$. The same is true using
(\ref{eq:gwmain}) as the definition of $C_{\mu\nu}^{\lambda,d}$.
Bertram, Ciocan-Fontaine, and Fulton \cite{b-c-f} (see also
\cite[{\S}3]{post}) give a formula for $C_{\mu\nu}^{\lambda,d}$ as an
alternating sum of Littlewood-Richardson coefficients, but again no
direct proof of positivity is known.
On the other hand, if we take $s_{\lambda/d/\mu}$ in more than $k$
variables, then it is still true (and not difficult to prove) that
$s_{\lambda/d/\mu}$ remains a symmetric function, but it need not be
Schur positive.

Theorem~\ref{thm:post} was used by Postnikov to obtain many properties
of Gromov-Witten invariants, some already known and some new.
For example, he gives a transparent explanation of a ``hidden''
symmetry of Gromov-Witten invariants, and he solves the previously
open problem of describing which powers $q^d$ appear with nonzero
coefficients in a quantum product $\sigma_\lambda\ast \sigma_\mu$ of
quantum Schubert classes. 
This latter problem is equivalent to determining which
$\mu^\vee/d/\lambda$ form a valid toric shape, where $\mu^\vee$ is the
complementary partition $(n-k-\mu_k,\dots, n-k-\mu_1)$.
For further details, see \cite{post}.

\section{Toric $h$-vectors and intersection cohomology.}
\label{sec:hvector} 
A \textbf{convex polytope} is the convex hull of finitely many points in
a Euclidean space. 
The subject of convex polytopes has undergone spectacular progress in
recent years.
New tools from commutative algebra, exterior algebra, algebraic
geometry, and other fields have led to solutions of previously
intractable problems.
Convex polytopes have also arisen in unexpected new areas and have
been applied to problems in these areas.
Computer science has raised a host of new problems related to convex
polytopes and has greatly increased our ability to handle specific
examples.
A good general reference is \cite{ziegler}.
We will discuss a recent breakthrough related to combinatorial
properties of convex polytopes, followed by a brief description of a
recent result on matroid complexes. 
Both results have a similar ``hard Lefschetz'' flavor.

Let $\cp$ be a $d$-dimensional (convex) polytope, or $d$-polytope for
short. 
Let $f_i = f_i(\cp)$ denote the number of $i$-dimensional faces of
$\cp$, where we set $f_{-1}=1$ (regarding the empty face as having
dimension $-1$). The vector $f(\cp)=(f_0, f_1,\dots, f_{d-1})$ is
called the $\bm{f}$\textbf{-vector} of ${\cal P}$.
The general problem of characterizing $f$-vectors of polytopes seem
hopeless, but when $\cp$ is \textbf{simplicial} (i.e., every proper
face of $\cp$ is a simplex) then a complete characterization is known
(the $\bm{g}$\textbf{-theorem} for simplicial polytopes).
Define the $\bm{h}$\textbf{-vector} $h(\cp)=(h_0,h_1,\dots,h_d)$ of a
simplicial polytope $\cp$ with $f$-vector $(f_0,\dots,f_{d-1})$ by
\beq 
 \sum_{i=0}^d f_{i-1}(x-1)^{d-i} = \sum_{i=0}^d h_i x^{d-i}. 
    \label{eq:hvector} 
\eeq
The $h$-vector and $f$-vector convey the same information. 
A vector $(h_0,h_1,\dots,h_d)\in \zz^{d+1}$ is the $h$-vector of a
simplicial $d$-polytope if and only if the following conditions are
satisfied:
  \be \item[(G$_1$)] $h_0=1$
  \item[(G$_2$)] $h_i=h_{d-i}$ for $0\leq i\leq d$ (the Dehn-Sommerville
   equations) 
  \item[(G$_3$)] $h_0\leq h_1\leq \cdots\leq h_{\lfloor d/2\rfloor}$ (the
    Generalized Lower Bound Conjecture, or GLBC)
  \item[(G$_4$)] Certain non-polynomial inequalities (which we call the
  $\bm{g}$\textbf{-in\-e\-qual\-i\-ties}) asserting that the
    differences $g_i:=h_i-h_{i-1}$ cannot grow too rapidly for $1\leq
 i\leq \lfloor d/2\rfloor$. We will not state these conditions here.
 \ee
The sufficiency of the above conditions was proved by Billera and Lee,
and the necessity by Stanley. 
A brief exposition of this result, with further references, appears in
\cite[{\S}3.1]{rs:cc}.
The basic idea of the proof of necessity is the following. First,
condition G$_1$ is trivial, and G$_2$ is well-known and not difficult
to prove.
To prove G$_3$ and G$_4$, slightly perturb the vertices of $\cp$ so
that they have rational coordinates.
Since $\cp$ is simplicial, small perturbations of the vertices do not
change the combinatorial type, and hence leave the $f$-vector and
$h$-vector invariant.
Once $\cp$ has rational vertices, we can construct a projective
algebraic variety $X_\cp$, the \textbf{toric variety} corresponding to
the normal fan $\Sigma_\cp$ of $\cp$, whose cohomology ring (say over
$\rr$) has the form
$$ 
H^*(X_\cp;\rr)=H^0(X_\cp;\rr)\oplus H^2(X_\cp;\rr)\oplus
\cdots\oplus H^{2d}(X_\cp;\rr), 
$$
where $\dim H^{2i}(X_\cp;\rr)=h_i(\cp)$. 
(The normal fan $\Sigma_\cp$ is defined below for any polytope $\cp$.)
Write for short $H^i=H^i(X_\cp;\rr)$.
Let $Y$ be a generic hyperplane section of $X_\cp$, with corresponding
cohomology class $\omega = [Y]\in H^2$.
The variety $X_\cp$ has sufficiently nice singularities (finite
quotient singularities) so that the hard Lefschetz theorem holds.
This means that if $i< d/2$, then the map
$\omega^{d-2i}:H^{2i}\rightarrow H^{2(d-i)}$ given by $u\mapsto
\omega^{d-2i} u$ is a bijection.
In particular, this implies that if $i<d/2$ then the map
$\omega:H^{2i}\rightarrow H^{2(i+1)}$ is injective.
Hence $\dim H^{2i}\leq \dim H^{2(i+1)}$ for $i<\lfloor d/2\rfloor$, so
G$_3$ follows.
To obtain G$_4$, we use the fact that $H^*(X_\cp;\rr)$ is a graded
$\rr$-algebra generated by $H^2$; details may be found in the
reference cited above.

What happens when we try to extend this reasoning to nonsimplicial
polytopes $\cp$? 
We can still define the variety $X_\cp$, but unless $\cp$ is
\textbf{rational} (i.e., has rational vertices) $X_\cp$ will not have
finite type (so will not be what is normally meant by an algebraic
variety) and little can be said.
Thus for now assume that $\cp$ has rational vertices.
Unfortunately $X_\cp$ now has more complicated singularities, and the
cohomology ring $H^*(X_\cp;\rr)$ is poorly behaved.
In particular, the Betti numbers $\dim H^i(X_\cp;\rr)$ depend on the
embedding of $\cp$ into the ambient Euclidean space.
However, the theory of \textbf{intersection cohomology} introduced by
Goresky and MacPherson \cite{g-m}\cite{g-m2} in 1980 yields ``nice''
cohomology for singular spaces.

In particular, if $\cp$ is any rational $d$-polytope then the
intersection cohomology of the toric variety $X_\cp$ has the form
$$ 
\ih(X_\cp;\rr) = \ih^0(X_\cp;\rr)\oplus
\ih^2(X_\cp;\rr)\oplus \cdots\oplus \ih^{2d}(X_\cp;\rr), 
$$
where each $\ih^{2i}(X_\cp;\rr)$ is a finite-dimensional real vector
space whose dimension $h_i$ depends only on the combinatorial type of
$\cp$ and not on its embedding into Euclidean space. 
The vector 
$$
h(\cp) = (h_0,h_1,\dots, h_d)
$$ 
is called the \textbf{toric} $\bm{h}$\textbf{-vector} (or
\textbf{generalized} $\bm{h}$\textbf{-vector}) of $\cp$.
If $\cp$ is simplicial, then intersection cohomology and singular
cohomology coincide, so the toric $h$-vector coincides with the
ordinary $h$-vector.

The combinatorial description of the toric $h$-vector is quite
subtle. 
For any polytope ${\cal P}$, define two polynomials $f(\cp,x)$ and
$g(\cp,x)$ recursively as follows.
If $\cp=\emptyset$ then $f(\emptyset,x) = g(\emptyset,x)=1$. If
$\dim\cp\geq 0$ then define
$$ 
f(\cp,x) = \sum_{{\cal Q}} g({\cal Q},x)(x-1)^{\dim\cp-\dim
     {\cal Q}-1}, 
$$
where ${\cal Q}$ ranges over all faces (including $\emptyset$) of $\cp$
except ${\cal Q}=\cp$. 
Finally if $\dim\cp =d\geq 0$ and $f(\cp,x) = h_0+h_1x+\cdots$, then
define
$$ 
g(\cp,x) = h_0+(h_1-h_0)x+(h_2-h_1)x^2+\cdots+(h_m-h_{m-1})
     x^m, 
$$
where $m=\lfloor d/2\rfloor$. 
It is easy to see that we have defined $f(\cp,x)$ and $g(\cp,x)$
recursively for all polytopes $\cp$.
For instance, let $\sigma_j$ denote a $j$-dimensional simplex and
${\cal C}_j$ a $j$-dimensional cube.
Suppose we have computed that $g(\sigma_0)=g(\sigma_1)=1$ and $g({\cal
C}_2,x)=1+x$.
Then
\beas 
 f({\cal C}_3,x) & = & 6(x+1)+12(x-1)+8(x-1)^2 +(x-1)^3\\
     & = & x^3+5x^2+5x+1, 
\eeas
and $g({\cal C}_3,x)=1+4x$.

It is easy to see that $\deg f(\cp,x)=d$, so $f(\cp,x)=h_0+h_1x+
\cdots +h_dx^d$. 
Then $h(\cp)=(h_0,h_1,\dots,h_d)$ is the toric
$h$-vector of $\cp$ when $\cp$ is rational, and we can \emph{define}
the toric $h$-vector of any polytope $\cp$ in this manner. 
As mentioned above, it coincides with the usual $h$-vector when $\cp$
is simplicial.
It is trivial that $h_0=1$, and it was first shown in
\cite[Thm.~2.4]{rs:genh} that $h_i=h_{d-i}$ for all $i$ (the
generalized Dehn-Sommerville equations).
Since $h_i=\dim\ih^{2i}(\xp;\rr)$ when $\cp$ is rational, we have
$h_i\geq 0$ in this case.
Moreover, $\ih(\xp;\rr)$ is a (graded) module over $H^*(\xp;\rr)$,
and the hard Lefschetz theorem continues to hold as follows: the map
$\omega:\ih^{2i}\rightarrow \ih^{2(i+1)}$ defined by $u\mapsto \omega
u$, where $\omega\in H^2$ as above, is injective for $i<d/2$.
Hence as before we get $h_0\leq h_1\leq \cdots\leq h_{\lfloor
d/2\rfloor}$.
Thus conditions G$_1$, G$_2$, and G$_3$ of the $g$-theorem continue to
hold for any rational polytope (using the toric $h$-vector).
However, intersection cohomology does not have a ring structure, and
it remains open whether G$_4$ holds for all rational polytopes.

Not all polytopes are rational, i.e., there exist polytopes $\cp$ for
which any polytope in $\rr^n$ with the same combinatorial type as
$\cp$ cannot have only rational vertices
\cite[Exam.~6.21]{ziegler}. 
It was conjectured in \cite{rs:genh} that conditions G$_3$--G$_4$ hold
for all convex polytopes.
An approach toward proving G$_3$ would be to construct a ``nice''
analogue of the toric variety $\xp$ when $\cp$ is nonrational.
It doesn't seem feasible to do this.
Instead we can try to construct an analogue of the intersection
cohomology $\ih(\xp;\rr)$ when $\cp$ is nonrational, i.e., a graded
$\rr$-algebra $\ih(\cp)=\ih^0(\cp)\oplus \ih^2(\cp)\oplus\cdots\oplus
\ih^{2d}(\cp)$ that becomes $\ih(\xp;\rr)$ when $\cp$ is rational and
which satisfies the two conditions:
  \be\item[(P$_1$)] $\dim \ih^{2i}(\cp)=h_i(\cp)$
  \item[(P$_2$)] $\ih(\cp)$ is a module over some ring containing an
element $l$ of degree 2 that satisfies the conditions of hard
Lefschetz theorem, i.e., for $i<d/2$ the map $l^{d-2i}:\ih^{2i}(\cp)
\rightarrow \ih^{2(d-i)}(\cp)$ is a bijection.
  \ee
\indent The first step in the above program was the definition of
$\ih(\cp)$ due to Barthel, Brasslet, Fieseler, and Kaup \cite{bbfk}
and to Bressler and Lunts \cite{b-l1}. 
The precise definition is rather technical; we include it here
(following \cite{karu}) so that even readers without the necessary
background will have some idea of its flavor.
With the polytope $\cp\subset \rr^n$ we can associate the
\textbf{normal fan} $\Sigma=\Sigma_\cp$, i.e., the set of all cones of
linear functions which are maximal on a fixed face of $\cp$ (e.g.,
\cite[Exam.~7.3]{ziegler}).
In other words, if $F$ is a nonempty face of ${\cal P}$ then define
$$ 
N_F := \left\{ f\in(\rr^n)^*\colon F\subseteq \{x\in\cp\colon
        f(x) =\max_{y\in\cp}f(y)\}\right\}. 
$$
The set of all cones $N_F$ forms the normal fan of $\cp$; it is a
\textbf{complete fan} in $(\rr^n)^*$, i.e., 
any two cones intersect in a common face of both, and the union of all
$N_F$'s is $(\rr^n)^*$. 
Define a sheaf ${\cal A}_\Sigma$, the \textbf{structure sheaf} of
$\Sigma$, as follows.
For each cone $\sigma\in\Sigma$ define the stalk ${\cal
A}_{\Sigma,\sigma}= \mathrm{Sym}(\mathrm{span}\,\sigma)^*$, the space
of polynomial functions on $\sigma$.
The restriction map ${\cal A}_{\Sigma,\sigma} \rightarrow {\cal
A}_\Sigma(\partial \sigma)$ is defined by restriction of functions
(where $\partial$ denotes boundary).
Thus ${\cal A}_\Sigma$ is a sheaf of algebras, naturally graded by
degree (where we define conewise linear functions to have degree 2),
so
$$
A_\Sigma = A_\Sigma^0\oplus A_\Sigma^2\oplus A_\Sigma^4\oplus\cdots.
$$
Let $A=\mathrm{Sym}\,(\rr^n)^*$ denote the space of polynomial
functions on all of $\rr^n$. 
Multiplication with elements of $A$ gives ${\cal A}_\Sigma$ the
structure of a sheaf of $A$-modules.

An \textbf{equivariant intersection cohomology sheaf} ${\cal L}_\Sigma$
of $\Sigma$ is a sheaf of ${\cal A}_\Sigma$-modules satisfying the
following three properties:
  \be\item[(E$_1$)] (normalization) ${\cal L}_{\Sigma,0}=\rr$
   \item[(E$_2$)] (local freeness) ${\cal L}_{\Sigma,\sigma}$ is a free
${\cal A}_{\Sigma,\sigma}$-module for any $\sigma\in\Sigma$.
   \item[(E$_3$)] (minimal flabbiness) Let $I$ be the ideal of $A$
generated by homogeneous linear functions, and for any $A$-module $M$
write $\overline{M}=M/IM$. 
Then modulo the ideal $I$ the restriction map induces an isomorphism
$$ \overline{{\cal L}}_{\Sigma,\sigma} \rightarrow \overline{{\cal
    L}_{\Sigma}(\partial \sigma)}. 
$$
  \ee

Equivariant intersection cohomology sheaves exist for any fan
$\Sigma$, and any two of them are isomorphic. 
Hence we may call a sheaf of ${\cal A}_\Sigma$-modules satisfying
E$_1$, E$_2$, and E$_3$ \emph{the} equivariant intersection homology
sheaf.
We have ${\cal L}_\Sigma \simeq {\cal A}_\Sigma$ if and only if the
fan $\Sigma$ is simplicial, in which case it coincides with the usual
equivariant cohomology sheaf.
(A good introduction to equivariant cohomology can be found in
\cite{brion}.)
The (non-equivariant) \textbf{intersection cohomology} of the fan
$\Sigma$ is defined to be the ${\cal A}$-module of global sections of
the intersection cohomology sheaf modulo the ideal $I$:
$$ 
\ih(\Sigma) = \overline{{\cal L}(\Sigma)}. 
$$
Since ${\cal L}(\Sigma)$ is a graded $A$-module and $I$ is a graded
ideal, $\ih(\Sigma)$ inherits a natural grading:
$$ 
\ih(\Sigma) = \ih(\Sigma)^0\oplus \ih(\Sigma)^2\oplus\cdots. 
$$
By definition ${\cal A}_\Sigma^2$ (the degree two part of ${\cal
A}_\Sigma$) consists of conewise linear functions on $\Sigma$.
The restriction $l_\sigma$ of $l$ to $\sigma$ is linear on $\sigma$
and hence extends to a unique linear function $l_\sigma\in A^2$.
A function $l\in{\cal A}_\Sigma^2$ is called \textbf{strictly convex}
if $l_\sigma(v)<l(v)$ for any maximal cone $\sigma\in\Sigma$ and any
$v\not\in\sigma$.
A complete fan $\Sigma$ is called \textbf{projective} if there exists
a strictly convex function $l\in {\cal A}_\Sigma^2$.

Now let $\cp$ be any convex polytope, and let $\Sigma_\cp$ denote the
normal fan of $\cp$. 
The normal fan $\Sigma_\cp$ is easy seen to be projective.
Write $\ih(\cp)$ for $\ih(\Sigma_\cp)$.

Bressler and Lunts \cite{b-l1} established a number of fundamental
properties of the intersection cohomology $\ih(\cp)$ of a convex
polytope $\cp$. 
They showed that if $\dim\cp=d$ then the grading of $\ih(\cp)$ has the
form
$$ 
\ih(\cp) = \ih^0(\cp)\oplus\ih^2(\cp)\oplus\cdots\oplus
        \ih^{2d}(\cp), 
$$
where each $\ih^{2i}(\cp)$ is a finite-dimensional vector space
satisfying Poincar\'e duality, so
$\ih^{2i}(\cp)\simeq \ih^{2(d-i)}(\cp)$. 
They conjectured that $\ih(\cp)$ has the hard Lefschetz property,
i.e., if $l\in{\cal A}_\Sigma^2$ is strictly convex, then for $i<d/2$
the map
$$ 
l^{d-i}:\ih^{2i}(\cp)\rightarrow\ih^{2(d-i)}(\cp) 
$$
is a bijection. 
They showed that if $\ih(\cp)$ does have the hard Lefschetz property,
then $\dim \ih^{2i}(\cp)=h_i(\cp)$.
Hence, as explained above, the toric $h$-vector of any polytope would
satisfy property G$_3$ (the GLBC).

The conjecture of Bressler and Lunts, that $\ih(\cp)$ has the hard
Lefschetz property, was first proved by Kalle Karu \cite{karu}.
Karu actually proves a stronger result, the
\textbf{Hodge-Riemann-Minkowski bilinear relations}.
To state this result, the Poincar\'e duality on intersection
cohomology gives a pairing
$$ 
\ih^{d-i}(\cp)\times \ih^{d+i}(\cp) \rightarrow\rr, 
$$
denoted $\langle x,y\rangle$. 
If $l\in {\cal A}_\Sigma^2$ is strictly convex, then we define a
quadratic form $Q_l$ on $\ih^{d-i}(\cp)$ by $Q_\ell(x)=\langle
l^ix,x\rangle$.
The \textbf{primitive intersection cohomology}
$\mathrm{IP}^{d-i}(\cp)$ is defined to be the kernel of the map
$$ 
l^{i+1}:\ih^{d-i}(\cp)\rightarrow\ih^{d+i+2}(\cp). 
$$
The Hodge-Riemann-Minkowski bilinear relations then state that the
quadratic form $(-1)^{(d-i)/2}Q_l$ is positive definite on
IP$^{d-i}(\cp)$ for all $i\geq 0$. 
(Note that IP$^{d-i}(\cp)=0$, or even $\ih^{d-i}(\cp)=0$, unless
$(d-i)/2$ is an integer.)
These relations were proved for simplicial polytopes by McMullen
\cite{mcm}.
A simpler proof was later given by Timorin \cite{tim}.
Very roughly, the idea of the proof of Karu is give a suitable
simplicial subdivision $\Delta$ of the fan $\Sigma_\cp$ and ``lift''
the Hodge-Riemann-Minkowski relations from $\Delta$ to $\Sigma$.

The proof of Karu has the defect that the Poincar\'e pairing $\langle
\cdot,\cdot\rangle$ depends on the choice of the subdivision $\Delta$
and on the embedding ${\cal L}(\Sigma)\subset {\cal L}(\Delta)$. 
In the words of Bressler and Lunts \cite{b-l2}, ``this ambiguity makes
the proof unnecessarily heavy and hard to follow.''
Bressler and Lunts define in this same paper a \emph{canonical}
pairing that considerably simplifies the proof of Karu.
They also show that Karu's pairing is independent of the choices made
to define it and in fact coincides with the canonical pairing.
Finally Barthel, Brasselet, Fieseler, and Kaup \cite{bbfk2} give a
``direct'' approach to the proof of Bressler and Lunts, replacing
derived categories with elementary sheaf theory and commutative
algebra.

We conclude this section by briefly mentioning a further recent result
concerning $h$-vectors and a variation of the hard Lefschetz theorem.
Undefined terminology may be found e.g.\ in \cite{rs:cc}.
A \textbf{matroid complex} is a (finite) abstract simplicial complex
$\Delta$ on a vertex set $V$ such that the restriction of $\Delta$ to
any subset of $V$ is pure, i.e., all maximal faces have the same
dimension.
The $f$-vector $f(\Delta)=(f_0,f_1,\dots,f_{d-1})$ is defined just as
for polytopes (where $\dim\Delta =d-1$), and the $h$-vector
$h(\Delta)=(h_0,h_1,\dots, h_d)$ is defined by equation
(\ref{eq:hvector}).
It is well known that matroid complexes are Cohen-Macaulay, so in
particular $h_i\geq 0$.
It is of considerable interest to obtain further conditions on
$h$-vectors of matroid complexes, and there are many results and
conjectures in this direction \cite[{\S}III.3]{rs:cc}.
In particular, Chari \cite[Cor.\ 4, part 2]{chari} has shown that if
$h_s\neq 0$ and $h_{s+1}=h_{s+2}=\cdots=h_d=0$, then $h_0\leq h_1\leq
\cdots\leq h_{\lfloor s/2\rfloor}$ and $h_i\leq h_{s-i}$ for $0\leq
i\leq \lfloor s/2\rfloor$.
Because matroid complexes are Cohen-Macaulay, there is a natural
graded $\rr$-algebra $A=A_0\oplus A_1\oplus\cdots\oplus A_s$ such that
$\dim A_i=h_i$ (namely, $A$ is the face ring of $\Delta$ modulo a
linear system of parameters).
This suggests that a ``half hard Lefschetz theorem'' might hold for
matroid complexes, viz., there is an element $l\in A_1$ such that the
map $l^{d-2i}:A_i\rightarrow A_{d-i}$ is \emph{injective} (rather than
bijective) for $0\leq i\leq \lfloor s/2\rfloor$.
In particular, this would imply the $g$-inequalities G$_4$ for the
sequence $h_0,h_1,\dots, h_{\lfloor s/2\rfloor}$.
This half hard Lefschetz theorem (for a generic linear system of
parameters) was proved by Ed Swartz \cite{swartz}.

\end{document}